\documentclass{article}

\usepackage{myarxiv}
\usepackage[utf8]{inputenc} % allow utf-8 input
\usepackage[T1]{fontenc}    % use 8-bit T1 fonts
\usepackage{hyperref}       % hyperlinks
\usepackage{url}            % simple URL typesetting
\usepackage{booktabs}       % professional-quality tables
\usepackage{amsfonts}       % blackboard math symbols
\usepackage{nicefrac}       % compact symbols for 1/2, etc.
\usepackage{microtype}      % microtypography
\usepackage{lipsum}

\usepackage{amssymb,amsmath,amsfonts,amsthm,enumerate}
\usepackage{subfigure}
\usepackage{graphicx}
\usepackage{epsfig}
\usepackage{float}
\usepackage[usenames]{color}
\usepackage{cite}
\usepackage{bm}

\newtheorem{theorem}{Theorem}[section]
\newtheorem{lemma}[theorem]{Lemma}%[section]
\newtheorem{definition}[theorem]{Definition}%[section]
\newtheorem{example}{Example}%[section]
\newtheorem{problem}[theorem]{Problem}%[section]
\numberwithin{equation}{section}

\def\BR{\mathbb R}

\def\cL{\mathcal L}
\def\cP{\mathcal P}

\def\T{\mathrm T}

\def\rd{\mathrm d}

\def\rdiv{\mathrm{div}}
\def\e{\mathrm e}

\def\De{\Delta}

\def\La{\Lambda}

\def\Om{\Omega}
\def\al{\alpha}

\def\ve{\varepsilon}
\def\te{\theta}
\def\ze{\zeta}
\def\ka{\kappa}
\def\la{\lambda}
\def\si{\sigma}
\def\vp{\varphi}
\def\om{\omega}
\def\f{\frac}
\def\nb{\nabla}
\def\ov{\overline}
\def\pa{\partial}

\def\wt{\widetilde}

\title{Unique Continuation Property with Partial Information for Two-Dimensional Anisotropic Elasticity Systems}

\author{
  Jin \uppercase{Cheng}\\
  School of Mathematical Sciences, Fudan University, Shanghai 200433, China;\\
  School of Mathematics, Shanghai University of Finance and Economics, Shanghai 200433, China.\\
  \texttt{jcheng@fudan.edu.cn}\\
   \And
  Yikan \uppercase{Liu}\thanks{Corresponding author.}\\
  Research Institute for Electronic Science, Hokkaido University, Sapporo 060-0812, Japan.\\
  \texttt{ykliu@es.hokudai.ac.jp}\\
  \AND
  Yanbo \uppercase{Wang}\\
  School of Mathematics, Shanghai University of Finance and Economics, Shanghai 200433, China.\\
  \texttt{yanbo68@hotmail.com} \\
   \And
  Masahiro \uppercase{Yamamoto}\\
  Graduate School of Mathematical Sciences, The University of Tokyo, Tokyo 153-8914, Japan;\\
  Academy of Romanian Scientists, Bucharest 050094, Romania;\\
  Peoples' Friendship University of Russia (RUDN University), Moscow 117198, Russian.\\
  \texttt{myama@ms.u-tokyo.ac.jp} \\
}

\begin{document}
\maketitle

\begin{abstract}
In this paper, we establish a novel unique continuation property for two-dimensional anisotropic elasticity systems with partial information. More precisely, given a homogeneous elasticity system in a connected open bounded domain, we investigate the unique continuation by assuming only the vanishing of one component of the solution in a subdomain. Using the corresponding Riemann function, we prove that the solution vanishes in the whole domain provided that the other component vanishes at one point up to its second derivatives. Further, we construct several examples showing the possibility of further reducing the additional information of the other component. This result possesses remarkable significance in both theoretical and practical aspects because the required data are almost halved for the unique determination of the whole solution.
\end{abstract}

\MRsubject{74B05\and 35B60}

\section{Introduction and Motivation}\label{sec-intro}

The unique continuation for partial differential equations (PDEs) usually refers to such properties that solutions to PDEs in small domains can uniquely determine those in larger ones. For a linear partial differential operator $\cP$ defined in a domain $\Om$, the unique continuation is concerned with the following type of problems.

\begin{problem}
If a solution $u$ to the homogeneous problem
\[
\cP u = 0\quad\mbox{in }\Om
\]
vanishes in some subdomain $\om\subset\Om$, does $u$ vanish in a larger subdomain of $\Om$?
\end{problem}

As one of the most fundamental properties of PDEs, the unique continuation characterizes the essential features of the equations. For instance, it is well known that Laplace and heat equations possess strong unique continuation in such a sense that the vanishing of a solution in an arbitrary subdomain $\om\subset\Om$ indicates its vanishing in the whole domain $\Om$. In contrast, due to the finite propagation speed of waves, the vanishing of a solution to a wave equation only implies its vanishing in a slightly larger ``cone-like'' domain. Regarding the classical theory on the unique continuation, we refer to the monographs \cite{H85,E86,R91} and the review \cite{SS87}.

Other than the qualitative analysis of PDEs, unique continuation properties (UCPs) are also frequently applied to such topics in PDEs as the determination of global solutions or other partial information by local measurements of solutions. Such applications are not only important theoretically, but also of practical significance in engineering and industry. For example, for the uniqueness of linear inverse problems, one can usually employ UCPs to extend the observation data to larger (sub) domains and then, in combination with some other methods, obtain the uniqueness of some coefficients in the PDEs (see e.g.\! \cite{CRV01,V08,Y09}). On the other hand, UCPs are also closely related to the controllability theory (see e.g.\! \cite{R78,I93,T95}).

Due to the limitation of observable information in practices, it remains a persistent and challenging issue in applied mathematics that how to effectively exploit the local information of solutions or equivalently, the maximal reduction of required measurements. Especially, owing to the mutual effects and influence of the components, a natural question arising in the treatment for coupling systems of PDEs is the possibility of acquiring the information of the whole systems by using only the local observation data of a part of components. In the context of inverse problems for PDEs, this is equivalent to the simultaneous identification of several unknown coefficients by partial information, and is related to controlling the whole system via less control force. On this direction, pioneering works have been carried out e.g.\! in \cite{CGR06,BCGY09,BCGT10,GT10,ABGT11,CGRY12,RC12,M13,BCGT14,DL16} for parabolic systems, \cite{CG10,DY19} for Schr\"odinger systems, and \cite{IIY03} for dynamical Lam\'e systems. It is witnessed in the above literature that the stronger the systems are coupled, the tighter the components are related and thus, the higher possible to ``visualize the whole leopard by just looking at one spot on it''.

The main purpose of this article is to establish a new UCP for two-dimensional elasticity systems with partial information. To articulate the problem, we formulate the system under consideration as follows. Let $\Om\subset\BR^2$ be a connected open bounded domain, and $\om(x_0,y_0)\subset\Om$ be an open neighborhood of some fixed point $(x_0,y_0)\in\Om$. We consider the following two-dimensional general anisotropic elasticity system for the displacement vector $\bm u=(u_1,u_2)^\T=(u_1(x,y),u_2(x,y))^\T$:
\begin{equation}\label{eq-gov-u}
\sum_{j,k,\ell=1}^2a_{ijk\ell}(x,y)\pa_j\pa_\ell u_k+\sum_{j,k=1}^2b_{ijk}(x,y)\pa_ku_j+\sum_{j=1}^2c_{ij}(x,y)u_j=0,\quad i=1,2\quad\mbox{in }\Om,
\end{equation}
where $(a_{ijk\ell})_{1\le i,j,k,\ell\le2}$ is a fully symmetric fourth order tensor, that is,
\begin{equation}\label{eq-sym}
a_{ijk\ell}=a_{jik\ell}=a_{k\ell ij},\quad\forall\,i,j,k,\ell=1,2\quad\mbox{on }\ov\Om\,.
\end{equation}
Here $(\,\cdot\,)^\T$ denotes the transpose of a matrix, and it is understood that $\pa_1=\pa_x=\f\pa{\pa x}$ and $\pa_2=\pa_y=\f\pa{\pa y}$. Various assumptions concerning the involved coefficients, especially $(a_{ijk\ell})$ will be specified later in Section \ref{sec-main}.

In this paper, we investigate the following problem.

\begin{problem}\label{prob-UCP}
Let $\bm u=(u_1,u_2)^\T$ satisfy the homogeneous elasticity system \eqref{eq-gov-u} in a connected open bounded domain $\Om\subset\BR^2,$ and $\om(x_0,y_0)\subset\Om$ be an open neighborhood of some fixed point $(x_0,y_0)\in\Om$. If the first component $u_1=0$ in $\om(x_0,y_0),$ does it imply $\bm u\equiv\bm0$ in $\Om$ under certain conditions$?$
\end{problem}

The researches on UCPs of elasticity systems can be traced back to about 25 years ago. For isotropic medium, Dehman and Robbiano \cite{DR93} first proved the UCP for Lam\'e systems with $C^\infty$ coefficients, and the coefficient regularity was reduced later in Ang et al.\! \cite{AITY98}. For anisotropic medium, Nakamura and Wang \cite{NW03} proved the UCP for elasticity systems with residual stress, and the same results were established later in \cite{NW06} with coefficients of Lipschitz continuity. For the corresponding strong UCP, we refer to \cite{AM01,L04,E06,LNW10,LNUW11}. It turns out that most existing literature focused on the reduction of coefficient regularity to generalize the results, whereas all of them required the vanishing of both components in a subdomain, i.e., $\bm u=(u_1,u_2)^\T=0$ in $\om(x_0,y_0)$. On the contrary, in this paper we only assume the vanishing of one component $u_1$ in $\om(x_0,y_0)$, and then take advantage of the coupling property of the system to transfer such information to the other component $u_2$, so that we can further deduce the vanishing of the solution $\bm u$ in the whole domain $\Om$. To the authors' best knowledge, such an aspect has remarkable novelty in the researches of UCPs of elasticity system, which also has potential real-world applications.

On this direction, only \cite{WLC17} treated the same problem for the isotropic prototype of \eqref{eq-gov-u}, that is, the Lam\'e system
\begin{equation}\label{eq-Lame}
\rdiv(\mu(\nb\bm u+(\nb\bm u)^\T))+\nb(\la\,\rdiv\,\bm u)=\bm0\quad\mbox{in }\Om,
\end{equation}
where $\mu=\mu(x,y)$ and $\la=\la(x,y)$ are the Lam\'e coefficients satisfying the strong elliptic condition \eqref{eq-cond-ell}, i.e., $\mu>0$ and $2\mu+\la>0$ on $\ov\Om$\,. One of the motivations of this paper is to generalize the result obtained in \cite{WLC17} to the anisotropic elasticity system \eqref{eq-gov-u}.

Theoretically, we are mainly inspired by the simplest situation of \eqref{eq-Lame} that both $\mu$ and $\la$ are constants. In this case, direct calculations and the vanishing assumption in Problem \ref{prob-UCP} immediately yield
\begin{equation}\label{eq-Lame-u2}
\begin{cases}
(\mu+\la)\pa_x\pa_y u_2=0,\\
\mu\,\pa_x^2u_2+(2\mu+\la)\pa_y^2u_2=0
\end{cases}\mbox{in }\om(x_0,y_0).
\end{equation}
Therefore, it is readily seen that as long as $\mu+\la\ne0$, the second component $u_2$ should take the form of
\begin{equation}\label{eq-Lame-sol}
u_2(x,y)=C_0+C_1x+C_2y+C_3\left(x^2-\f\mu{2\mu+\la}y^2\right)\quad\mbox{in }\om(x_0,y_0),
\end{equation}
where $C_0,C_1,C_2,C_3$ are constants determined by extra information e.g.\! on $u_2$ and its derivatives at $(x_0,y_0)$. Hence, if such information implies $C_0=C_1=C_2=C_3=0$ and thus $u_2=0$ in $\om(x_0,y_0)$, then a direct application of some existing UCPs such as Lemma \ref{lem-NW} gives a positive answer to Problem \ref{prob-UCP}.

Regardless of its simplicity, the isotropic Lam\'e system \eqref{eq-Lame} with constant coefficients provides informative messages to our problem. First, the vanishing of one component is indeed a rather strong restriction in such a sense that it forces the other to satisfy an overdetermined system like \eqref{eq-Lame-u2}. Further, the equations in \eqref{eq-Lame-u2} are in principle of hyperbolic and elliptic types, which is expected to be inherited in the general anisotropic case \eqref{eq-gov-u}. Second, the overdetermined system \eqref{eq-Lame-u2} is still not enough to determine the solution, but the required additional information should be definitely much less than that in traditional UCPs.

The rest of this article is organized as follows. In Section \ref{sec-main}, we recall a representative UCP of elasticity systems and state the main result (Theorem \ref{thm-UCP}) answering Problem \ref{prob-UCP} along with remarks and examples. Then Section \ref{sec-proof} is devoted to the proof of Theorem \ref{thm-UCP}, and we provide further explanations and construct several examples on reducing the additional information in Section \ref{sec-egrem}. Finally, we close this paper with concluding remarks in Section \ref{sec-con}.

\section{Preliminary and the Main Result}\label{sec-main}

We start with some general settings in Problem \ref{prob-UCP}. Throughout this paper, we restrict ourselves to considering only classical solutions $\bm u\in(C^2(\ov\Om))^2$ satisfying the governing system \eqref{eq-gov-u} pointwisely, which is guaranteed by suitably chosen coefficients and boundary conditions. As the well-posedness issue of \eqref{eq-gov-u} is not the main focus of this paper, we omit the details and refer the interested readers to e.g.\! Li and Qin \cite{LQ12}.

Next, we fix the assumptions on the coefficients that appear in the governing system \eqref{eq-gov-u}. Besides the condition \eqref{eq-sym} for symmetry, we recall two kinds of conditions on the tensor $(a_{ijk\ell})$ for the well-posedness of \eqref{eq-gov-u} appearing frequently in the literature (see e.g. \cite{LQ12}).

\begin{definition}\label{def-tensor}
Let $(a_{ijk\ell})$ be a fully symmetric tensor satisfying \eqref{eq-sym}.

{\rm(a)}\ \ We say that $(a_{ijk\ell})$ satisfies the strong elliptic condition, if there exists a constant $\ka_0>0$ such that for all $(x,y)\in\ov\Om$ and all $(\xi_1,\xi_2),(\eta_1,\eta_2)\in\BR^2$, there holds
\begin{equation}\label{eq-cond-ell}
\sum_{i,j,k,\ell=1}^2a_{ijk\ell}(x,y)\xi_i\eta_j\xi_k\eta_\ell\ge\ka_0(\xi_1^2+\xi_2^2)(\eta_1^2+\eta_2^2).
\end{equation}

{\rm(b)}\ \ We say that $(a_{ijk\ell})$ satisfies the strong convexity condition, if there exists a constant $\ka_1>0$ such that for all $(x,y)\in\ov\Om$ and all symmetric matrices $(e_{ij})_{1\le i,j\le 2}$, there holds
\begin{equation}\label{eq-cond-cvx}
\sum_{i,j,k,\ell=1}^2a_{ijk\ell}(x,y)e_{ij}e_{k\ell}\ge\ka_1\sum_{i,j=1}^2e_{ij}^2.
\end{equation}
\end{definition}

Obviously, the strong convexity condition is stronger than the strong elliptic condition. Throughout this paper, we basically assume the strong elliptic condition \eqref{eq-cond-ell} unless otherwise specified.

Concerning the regularity, we assume that for all $i,j,k,\ell=1,2$, the second order coefficients $a_{ijk\ell}$ are locally Lipschitz, and the lower order ones $b_{ijk},c_{ij}$ are locally bounded in $\Om$. In accordance with the notations in \cite{NW06}, we can rewrite \eqref{eq-gov-u} as
\begin{equation}\label{eq-gov-u'}
\La_{11}\pa_x^2\bm u+\La_{12}\pa_x\pa_y\bm u+\La_{22}\pa_y^2\bm u+D(\bm u)=\bm0\quad\mbox{in }\Om,
\end{equation}
where
\[
\La_{11}=\begin{pmatrix}
a_{1111} & a_{1112}\\
a_{1112} & a_{1212}
\end{pmatrix},\quad\La_{12}=\begin{pmatrix}
2a_{1112} & a_{1212}+a_{1122}\\
a_{1212}+a_{1122} & 2a_{1222}
\end{pmatrix},\quad\La_{22}=\begin{pmatrix}
a_{1212} & a_{1222}\\
a_{1222} & a_{2222}
\end{pmatrix},
\]
and $D$ is a first order differential operator with locally bounded coefficients.

As a representative of existing UCPs for two-dimensional anisotropic elasticity systems, we recall the following result.

\begin{lemma}\label{lem-NW}
Let $\bm u$ satisfy \eqref{eq-gov-u}, $\om(x,y)\subset\Om$ be a sufficiently small open neighborhood of $(x,y)\in\Om,$ and $(\te,(z_1,z_2)^\T)$ be a local eigenpair of the quadratic pencil $\La_{11}\ze^2+\La_{12}\ze+\La_{22},$ i.e.
\[
(\La_{11}\te^2+\La_{12}\te+\La_{22})(z_1,z_2)^\T=\bm0\quad\mbox{in }\om(x,y).
\]
Assume that for all $(x,y)\in\Om,$ the functions $\te,(z_1,z_2)^\T$ are Lipschitz and the matrix function $\begin{pmatrix}
z_1 & \ov{z_1}\\
z_2 & \ov{z_2}
\end{pmatrix}$ is nonsingular in $\om(x,y)$. Then $\bm u=(u_1,u_2)^\T=\bm0$ in $\om(x_0,y_0)$ implies $\bm u\equiv\bm0$ in $\Om$.
\end{lemma}

The conclusion in the above lemma is a straightforward corollary of Nakamura and Wang \cite[Theorem 2.4]{NW06}, which only required $\bm u$ to vanish at some $(x_0,y_0)\in\Om$ of any exponential order. On the other hand, note that \cite{NW06} only treated the divergence form
\[
\nb\cdot((a_{ijk\ell}(x,y))\nb\bm u)=\bm0\quad\mbox{in }\Om
\]
without any other lower order terms. However, the argument there works in the more general formulation \eqref{eq-gov-u} because the lower order terms collected in $D(\bm u)$ are unimportant.

With the above preparation, keeping the motivation in Section \ref{sec-intro} in mind, now we state the main result of this article.

\begin{theorem}\label{thm-UCP}
Let $(x_0,y_0)\in\Om$ be arbitrarily fixed. Let $\bm u=(u_1,u_2)^\T\in(C^2(\ov\Om))^2$ be a classical solution of $\eqref{eq-gov-u},$ and $\om(x_0,y_0)\subset\Om$ be an open neighborhood of $(x_0,y_0)$. Assume that the fully symmetric tensor $(a_{ijk\ell})$ satisfies the strong elliptic condition $\eqref{eq-cond-ell},$ and the assumptions of Lemma $\ref{lem-NW}$ are fulfilled throughout $\Om$. If
\begin{equation}\label{eq-cond-hyper}
\De:=(a_{1212}+a_{1122})^2-4a_{1112}a_{1222}>0\quad\mbox{on }\ov{\om(x_0,y_0)}\,,
\end{equation}
then $u_1=0$ in $\om(x_0,y_0)$ implies $\bm u\equiv\bm0$ in $\Om$ with the additional information
\begin{equation}\label{eq-vanish-u2}
u_2=\pa_x u_2=\pa_y u_2=\pa_x^2u_2=\pa_y^2u_2=0\quad\mbox{at }(x_0,y_0).
\end{equation}
\end{theorem}

The new UCP in the above theorem gives an affirmative answer to Problem \ref{prob-UCP}, which enables us to reduce almost half information of $\bm u$ under certain assumptions. In replace of taking efforts to observe $u_2$ throughout $\om(x_0,y_0)$, it suffices to perform single point observation such as \eqref{eq-vanish-u2}. Nevertheless, the essential assumption \eqref{eq-cond-hyper} looks unusual in the related literature, and thus some explanations are necessary to demonstrate its feasibility.

In fact, the propose of \eqref{eq-cond-hyper} is simply to preserve the hyperbolic-elliptic structure as the simplest case \eqref{eq-Lame-u2}, which seems natural and reasonable. On the other hand, in view of the formulation \eqref{eq-gov-u'}, we observe that \eqref{eq-cond-hyper} is equivalent to assuming
\[
\det\La_{12}=-\De=4a_{1112}a_{1222}-(a_{1212}+a_{1122})^2<0\quad\mbox{on }\ov{\om(x_0,y_0)}\,,
\]
which is, of course, absent in the assumption of Lemma \ref{lem-NW}. Since such an extra assumption is only imposed in the subset $\ov{\om(x_0,y_0)}\,$, we understand that \eqref{eq-cond-hyper} is basically independent of the assumption of Lemma \ref{lem-NW} in general.

Next, we further compare \eqref{eq-cond-hyper} with the strong convexity condition \eqref{eq-cond-cvx} in two special cases. In the isotropic medium, we have
\[
a_{1111}=a_{2222}=2\mu+\la,\quad a_{1212}=\mu,\quad a_{1122}=\la,\quad a_{1112}=a_{1222}=0.
\]
Then it is readily seen that condition \eqref{eq-cond-hyper} reduces to
\begin{equation}\label{eq-cond-Lame}
\mu+\la\ne0\quad\mbox{on }\ov{\om(x_0,y_0)}\,.
\end{equation}
Meanwhile, it is easily verified that the condition \eqref{eq-cond-cvx} in the isotropic case reads
\[
\mu>0,\ \mu+\la>0\quad\mbox{on }\ov\Om\,,
\]
which is definitely stronger than \eqref{eq-cond-Lame}. Indeed, if $\mu+\la=0$ in $\om(x_0,y_0)$, then the Lam\'e system \eqref{eq-Lame} becomes
\begin{equation}\label{eq-Lame-weak}
\begin{cases}
\rdiv(\mu\nb u_1)+(\pa_y\mu)\pa_x u_2-(\pa_x\mu)\pa_y u_2=0,\\
(\pa_x\mu)\pa_y u_1-(\pa_y\mu)\pa_x u_1+\rdiv(\mu\nb u_2)=0
\end{cases}\mbox{in }\om(x_0,y_0).
\end{equation}
Even if $u_1=0$ in $\om(x_0,y_0)$, the remaining equations for $u_2$ fail to indicate $u_2=0$ in $\om(x_0,y_0)$ with the additional information \eqref{eq-vanish-u2}. Because now \eqref{eq-Lame-weak} is at most a weakly coupling system, it is understood that the connection between $u_1$ and $u_2$ is not strong enough for $u_1$ to determine $u_2$ uniquely. By the way, since the assumption in Lemma \ref{lem-NW} holds automatically in the isotropic case (see \cite[Example 2]{NW06}), we conclude that \eqref{eq-cond-Lame} (or \eqref{eq-cond-hyper}) is an intermediate condition between \eqref{eq-cond-ell} and \eqref{eq-cond-cvx} in this case.

Now we turn to the orthotropic medium, that is, only $a_{1112}=a_{1222}=0$. In this case, the condition \eqref{eq-cond-hyper} becomes
\begin{equation}\label{eq-cond-ortho}
a_{1212}+a_{1122}\ne0\quad\mbox{on }\ov{\om(x_0,y_0)}\,,
\end{equation}
which preserves the form of \eqref{eq-cond-Lame} as that in the isotropic case. However, the strong convexity condition \eqref{eq-cond-cvx} does not necessarily imply \eqref{eq-cond-hyper} any more; a simple counterexample can be
\[
a_{1111}=2,\quad a_{1212}=a_{2222}=1,\quad a_{1122}=-1.
\]
Therefore, it can be inferred that \eqref{eq-cond-hyper} is also basically independent of \eqref{eq-cond-cvx} in general.

Next we discuss the additional data \eqref{eq-vanish-u2}, which consists of $5$ pieces of information in the single point observation of $u_2$. Recall that in the simplest case of constant Lam\'e coefficients, $u_2$ takes the form of \eqref{eq-Lame-sol}, which is uniquely determined by $4$ constants. Then it is readily seen that the extra information on $u_2$ can be reduced to $4$ pieces as
\begin{equation}\label{eq-info-u2-4}
u_2=\pa_x u_2=\pa_y u_2=0,\quad\pa_x^2u_2=0\mbox{ or }\pa_y^2u_2=0\quad\mbox{at }(x_0,y_0).
\end{equation}
Actually, it follows immediately from \eqref{eq-Lame-u2} that either $\pa_x^2u_2=0$ or $\pa_y^2u_2=0$ implies the vanishing of the other at $(x_0,y_0)$. However, we will provide examples in Section \ref{sec-egrem} showing that $4$ pieces of information such as \eqref{eq-info-u2-4} are possibly insufficient in general. On the opposite side, we will also construct several examples with special choices of coefficients, such that the extra pieces of information on $u_2$ can be less than $4$. Consequently, it turns out that \eqref{eq-vanish-u2} is a sufficient condition to conclude Theorem \ref{thm-UCP}.

\section{Proof of Theorem \ref{thm-UCP}}\label{sec-proof}

Thanks to the powerful UCP result in Lemma \ref{lem-NW}, it is sufficient to prove $u_2=0$ in some neighborhood of $(x_0,y_0)$, which is not necessarily $\om(x_0,y_0)$. We divide the proof into two steps.\medskip

{\bf Step 1.}\ \ First we perform changes of variables to transform the equations for $u_2$ into specified forms. Simply substituting $u_1=0$ in $\om(x_0,y_0)$ into the governing system \eqref{eq-gov-u} yields
\begin{align}
& a_{1112}\pa_x^2u_2+(a_{1212}+a_{1122})\pa_x\pa_y u_2+a_{1222}\pa_y^2u_2+b_{121}\pa_x u_2+b_{122}\pa_y u_2+c_{12}u_2=0,\label{eq-u2-hyper}\\
& a_{1212}\pa_x^2u_2+2a_{1222}\pa_x\pa_y u_2+a_{2222}\pa_y^2u_2+b_{221}\pa_x u_2+b_{222}\pa_y u_2+c_{22}u_2=0\label{eq-u2-ell}
\end{align}
in $\om(x_0,y_0)$. Due to the essential assumption \eqref{eq-cond-hyper}, we see that equation \eqref{eq-u2-hyper} is of hyperbolic type. Therefore, there exist invertible changes of variables
\[
s=s(x,y),\quad t=t(x,y)
\]
such that upon introducing $w(s(x,y),t(x,y)):=u_2(x,y)$, the auxiliary function $w$ satisfies
\begin{equation}\label{eq-w-hyper}
\pa_s\pa_t w+B_{11}\pa_s w+B_{12}\pa_t w+C_1w=0
\end{equation}
on $\wt\om:=\{(s(x,y),t(x,y))\in\BR^2;\,(x,y)\in\ov{\om(x_0,y_0)}\}$. More precisely, we accomplish the changes of variables in the following way.
\begin{itemize}
\item If $a_{1112}=a_{1222}=0$, then changes of variables are not necessary, and we simply put
\[
B_{11}=\f{b_{121}}{a_{1212}+a_{1122}},\quad B_{12}=\f{b_{122}}{a_{1212}+a_{1122}},\quad C_1=\f{c_{12}}{a_{1212}+a_{1122}}.
\]
Note that in this case, the condition \eqref{eq-cond-hyper} guarantees $a_{1212}+a_{1122}\ne0$ (see \eqref{eq-cond-ortho}).
\item If $a_{1112}\ne0$, then we choose $\f{\pa_x s}{\pa_y s}$ and $\f{\pa_x t}{\pa_y t}$ as the two distinct solutions of the quadratic equation
\[
a_{1112}\eta^2+(a_{1212}+a_{1122})\eta+a_{1222}=0,
\]
that is,
\[
\f{\pa_x s}{\pa_y s}=-\f{a_{1212}+a_{1122}-\sqrt\De}{2a_{1112}},\quad\f{\pa_x t}{\pa_y t}=-\f{a_{1212}+a_{1122}+\sqrt\De}{2a_{1112}}.
\]
Here we recall $\De=(a_{1212}+a_{1122})^2-4a_{1112}a_{1222}>0$ on $\ov{\om(x_0,y_0)}$ (see \eqref{eq-cond-hyper}).
\item Similarly, if $a_{1222}\ne0$, then we choose $\f{\pa_y s}{\pa_x s}$ and $\f{\pa_y t}{\pa_x t}$ as the two distinct solutions of
\[
a_{1222}\eta^2+(a_{1212}+a_{1122})\eta+a_{1112}=0,
\]
that is,
\[
\f{\pa_y s}{\pa_x s}=-\f{a_{1212}+a_{1122}-\sqrt\De}{2a_{1222}},\quad\f{\pa_y t}{\pa_x t}=-\f{a_{1212}+a_{1122}+\sqrt\De}{2a_{1222}}.
\]
\end{itemize}
Recall that we assumed the locally Lipschitz continuity of all $a_{ijk\ell}$ in $\Om$ and the strict positivity of $\De$ on $\ov{\om(x_0,y_0)}\,$. Therefore, in all the above cases we can suitably specify $s(x,y)$ and $t(x,y)$, such that their first order derivatives are also locally Lipschitz and their second order ones are locally bounded on $\ov{\om(x_0,y_0)}\,$. Moreover, we can suppose $s(x_0,y_0)=t(x_0,y_0)=0$ so that the origin $(0,0)\in\wt\om$ without loss of generality, or otherwise it suffices to replace $s$ and $t$ by
\[
s(x,y)-s(x_0,y_0),\quad t(x,y)-t(x_0,y_0)
\]
respectively.

Now we turn to the second equation \eqref{eq-u2-ell}. Taking $\xi_1=0$, \ $\xi_2\ne0$ and $\eta_2\ne0$ in the strong elliptic condition \eqref{eq-cond-ell}, we obtain
\[
a_{1212}\eta^2+2a_{1222}\eta+a_{2222}>0\quad\mbox{on }\ov\Om,\quad\eta:=\f{\eta_1}{\eta_2}\in\BR.
\]
This implies $a_{1222}^2-a_{1212}a_{2222}<0$ on $\ov\Om$ and thus equation \eqref{eq-u2-ell} is automatically elliptic. In correspondence with the changes of variables, \eqref{eq-u2-ell} is transformed as
\begin{equation}\label{eq-w-ell}
\cL_{s,t}w:=A_{11}\pa_s^2w+2A_{12}\pa_s\pa_t w+A_{22}\pa_t^2w+B_{21}\pa_s w+B_{22}\pa_t w+C_2w=0.
\end{equation}
By the invertibility of $(s(x,y),t(x,y))$, we have
\begin{equation}\label{eq-invert}
\det\begin{pmatrix}
\pa_x s & \pa_x t\\
\pa_y s & \pa_y t
\end{pmatrix}\ne0\quad\mbox{on }\ov{\om(x_0,y_0)}\,,
\end{equation}
and equation \eqref{eq-w-ell} keeps the ellipticity because
\[
A_{12}^2-A_{11}A_{22}\propto\left(a_{1222}^2-a_{1212}a_{2222}\right)\left(\det\begin{pmatrix}
\pa_x s & \pa_x t\\
\pa_y s & \pa_y t
\end{pmatrix}\right)^2<0\quad\mbox{on }\wt\om.
\]

Now that $s(x,y)$ and $t(x,y)$ are well-defined up to their second order derivatives, direct calculations immediately give the relations
\begin{equation}\label{eq-u2-w-1}
\begin{aligned}
\pa_x u_2 & =(\pa_x s)\pa_s w+(\pa_x t)\pa_t w,\\
\pa_y u_2 & =(\pa_y s)\pa_s w+(\pa_y t)\pa_t w
\end{aligned}
\end{equation}
and further,
\begin{equation}\label{eq-u2-w-2}
\begin{aligned}
\pa_x^2u_2 & =(\pa_x s)^2\pa_s^2w+2(\pa_x s)(\pa_x t)\pa_s\pa_t w+(\pa_x t)^2\pa_t^2w+(\pa_x^2s)\pa_s w+(\pa_x^2t)\pa_t w,\\
\pa_x\pa_y u_2 & =(\pa_x s)(\pa_y s)\pa_s^2w+\{(\pa_x s)\pa_y t+(\pa_y s)\pa_x t\}\pa_s\pa_t w+(\pa_x t)(\pa_y t)\pa_t^2w\\
& \quad\,+(\pa_x\pa_y s)\pa_s w+(\pa_x\pa_y t)\pa_t w,\\
\pa_y^2u_2 & =(\pa_y s)^2\pa_s^2w+2(\pa_y s)(\pa_y t)\pa_s\pa_t w+(\pa_y t)^2\pa_t^2w+(\pa_y^2s)\pa_s w+(\pa_y^2t)\pa_t w.
\end{aligned}
\end{equation}
Based on the above relations, we can calculate all the coefficients involved in \eqref{eq-w-hyper} and \eqref{eq-w-ell}, and we omit the details here.

Now we convert the additional information \eqref{eq-vanish-u2} of $u_2$ at $(x_0,y_0)$ into that of $w$ at $(0,0)$. First, $w(0,0)=0$ is trivial. To see $\pa_s w=\pa_t w=0$ at $(0,0)$, it suffices to take $(x,y)=(x_0,y_0)$ in \eqref{eq-u2-w-1} and employ the invertibility condition \eqref{eq-invert}. For the second order derivatives, we shall plug \eqref{eq-vanish-u2} into the original equations \eqref{eq-u2-hyper}--\eqref{eq-u2-ell}, which gives
\[
(a_{1212}+a_{1122})\pa_x\pa_y u_2=a_{1222}\pa_x\pa_y u_2=0\quad\mbox{at }(x_0,y_0).
\]
According to condition \eqref{eq-cond-hyper}, it reveals that $a_{1212}+a_{1122}$ and $a_{1222}$ cannot vanish simultaneously at $(x_0,y_0)$, indicating $\pa_x\pa_y u_2(x_0,y_0)=0$. Now that
\[
\pa_x^2u_2=\pa_x\pa_y u_2=\pa_y^2u_2=0\mbox{ at }(x_0,y_0),\quad\pa_s w=\pa_t w=0\mbox{ at }(0,0),
\]
taking $(x,y)=(x_0,y_0)$ in \eqref{eq-u2-w-2} yields a homogeneous $3\times3$ linear system
\[
\begin{cases}
(\pa_x s)^2\pa_s^2w+2(\pa_x s)(\pa_x t)\pa_s\pa_t w+(\pa_x t)^2\pa_t^2w=0,\\
(\pa_x s)(\pa_y s)\pa_s^2w+\{(\pa_x s)\pa_y t+(\pa_y s)\pa_x t\}\pa_s\pa_t w+(\pa_x t)(\pa_y t)\pa_t^2w=0, & \mbox{at }(x_0,y_0).\\
(\pa_y s)^2\pa_s^2w+2(\pa_y s)(\pa_y t)\pa_s\pa_t w+(\pa_y t)^2\pa_t^2w=0
\end{cases}
\]
Employing the invertibility condition \eqref{eq-invert} again, we find
\[
\det\begin{pmatrix}
(\pa_x s)^2 & 2(\pa_x s)\pa_x t & (\pa_x t)^2\\
(\pa_x s)\pa_y s & (\pa_x s)\pa_y t+(\pa_y s)\pa_x t & (\pa_x t)\pa_y t\\
(\pa_y s)^2 & 2(\pa_y s)\pa_y t & (\pa_y t)^2
\end{pmatrix}=\left(\det\begin{pmatrix}
\pa_x s & \pa_x t\\
\pa_y s & \pa_y t
\end{pmatrix}\right)^3\ne0\quad\mbox{on }\ov{\om(x_0,y_0)}\,.
\]
Eventually, all derivatives of $w$ up to the second order vanish at $(0,0)$, i.e.,
\begin{equation}\label{eq-vanish-w}
w=\pa_s w=\pa_t w=\pa_s^2w=\pa_s\pa_t w=\pa_t^2w=0\quad\mbox{at }(0,0).
\end{equation}

{\bf Step 2.}\ \ We shall prove the claim by showing that there only exists a trivial solution of the overdetermined system \eqref{eq-w-hyper}--\eqref{eq-w-ell} with the additional information \eqref{eq-vanish-w}.

Let $\wt\om_1:=[-\ve,\ve]^2\subset\wt\om$ be a sufficiently small closed neighborhood square of $(0,0)$. According to the following key lemma, we can represent the solution of \eqref{eq-w-hyper} on $\wt\om_1$ by using the corresponding Riemann function.

\begin{lemma}[Vekua \cite{V68}]\label{lem-Riemann}
{\rm(a)}\ \ For any fixed $(\xi,\eta)\in\wt\om_1,$ there exists a unique Riemann function $R(s,t,\xi,\eta)$ of \eqref{eq-w-hyper} satisfying
\begin{equation}\label{eq-Riemann}
\begin{aligned}
&R(s,t,\xi,\eta)  -\int_\xi^sB_{12}(\si,t)R(\si,t,\xi,\eta)\,\rd\si-\int_\eta^tB_{11}(s,\tau)R(s,\tau,\xi,\eta)\,\rd\tau\\
& +\int_\xi^s\rd\si\int_\eta^tC_1(\si,\tau)R(\si,\tau,\xi,\eta)\,\rd\tau=1,\quad(s,t)\in\wt\om_1.
\end{aligned}
\end{equation}

{\rm(b)}\ \ The solution $w$ of \eqref{eq-w-hyper} takes the form of
\begin{equation}\label{eq-sol-w}
w(s,t)=w(0,0)R(0,0,s,t)+\int_0^s R(\si,0,s,t)\vp(\si)\,\rd\si+\int_0^t R(0,\tau,s,t)\psi(\tau)\,\rd\tau
\end{equation}
for $(s,t)\in\wt\om_1,$ where
\begin{equation}\label{eq-vp-psi}
\vp(s):=\pa_s w(s,0)+B_{12}(s,0)w(s,0),\quad
\psi(t):=\pa_t w(0,t)+B_{11}(0,t)w(0,t).
\end{equation}
\end{lemma}

Our strategy is to show $\vp=\psi=0$ on $[-\ve,\ve]$ and thus $w=0$ on $\wt\om_1$ by substituting the explicit solution \eqref{eq-sol-w} into the elliptic equation \eqref{eq-w-ell}. To this end, we first utilize the vanishing condition \eqref{eq-vanish-w} to obtain
\begin{equation}\label{eq-vp-psi-0}
\vp(0)=\vp'(0)=\psi(0)=\psi'(0)=0,
\end{equation}
and simplify \eqref{eq-sol-w} as
\[
w(s,t)=\int_0^s R(\si,0,s,t)\vp(\si)\,\rd\si+\int_0^t R(0,\tau,s,t)\psi(\tau)\,\rd\tau.
\]
Then direct calculations yield
\begin{align*}
\pa_s w(s,t) & =R(s,0,s,t)\vp(s)+\int_0^s\pa_s R(\si,0,s,t)\vp(\si)\,\rd\si+\int_0^t\pa_s R(0,\tau,s,t)\psi(\tau)\,\rd\tau,\\
\pa_t w(s,t) & =R(0,t,s,t)\psi(t)+\int_0^s\pa_t R(\si,0,s,t)\vp(\si)\,\rd\si+\int_0^t\pa_t R(0,\tau,s,t)\psi(\tau)\,\rd\tau,\\
\pa_s^2w(s,t) & =R(s,0,s,t)\vp'(s)+(\pa_s+2\pa_\xi)R(s,0,\xi,t)|_{\xi=s}\,\vp(s)\\
& \quad\,+\int_0^s\pa_s^2R(\si,0,s,t)\vp(\si)\,\rd\si+\int_0^t\pa_s^2R(0,\tau,s,t)\psi(\tau)\,\rd\tau,\\
\pa_s\pa_t w(s,t) & =\pa_t R(s,0,s,t)\vp(s)+\pa_s R(0,t,s,t)\psi(t)\\
& \quad\,+\int_0^s\pa_s\pa_t R(\si,0,s,t)\vp(\si)\,\rd\si+\int_0^t\pa_s\pa_t R(0,\tau,s,t)\psi(\tau)\,\rd\tau,\\
\pa_t^2w(s,t) & =R(0,t,s,t)\psi'(t)+(\pa_t+2\pa_\eta)R(0,t,s,\eta)|_{\eta=t}\,\psi(t)\\
& \quad\,+\int_0^s\pa_t^2R(\si,0,s,t)\vp(\si)\,\rd\si+\int_0^t\pa_t^2R(0,\tau,s,t)\psi(\tau)\,\rd\tau.
\end{align*}
Plugging the above representations into \eqref{eq-w-ell} and rearranging to separate the terms with $\vp$ and $\psi$, we arrive at
\begin{equation}\label{eq-ell-vp-psi}
\begin{aligned}
& A_{11}(s,t)R(s,0,s,t)\vp'(s)+P(s,t)\vp(s)+\int_0^s\cL_{s,t}R(\si,0,s,t)\vp(\si)\,\rd\si\\
& +A_{22}(s,t)R(0,t,s,t)\psi'(t)+Q(s,t)\psi(t)+\int_0^t\cL_{s,t}R(0,\tau,s,t)\psi(\tau)\,\rd\tau=0,
\end{aligned}
\end{equation}
where
\begin{align*}
P(s,t) & :=A_{11}(s,t)(\pa_s+2\pa_\xi)R(s,0,\xi,t)|_{\xi=s}+2A_{12}(s,t)\pa_t R(s,0,s,t)+B_{21}(s,t)R(s,0,s,t),\\
Q(s,t) & :=A_{22}(s,t)(\pa_t+2\pa_\eta)R(0,t,s,\eta)|_{\eta=t}+2A_{12}(s,t)\pa_s R(0,t,s,t)+B_{22}(s,t)R(0,t,s,t),
\end{align*}
and $\cL_{s,t}$ is the elliptic operator with respect to $s,t$ defined by \eqref{eq-w-ell}.

By the ellipticity of \eqref{eq-w-ell}, we see that $A_{11}A_{22}>0$ on $\wt\om$ and thus both $A_{11}$ and $A_{22}$ never vanish on $\wt\om$. In addition, taking $s=\xi$ and $t=\eta$ in \eqref{eq-Riemann}, we have $R(s,t,s,t)=1$ for all $(s,t)\in\wt\om_1$. Taking $t=0$ in \eqref{eq-ell-vp-psi} and using \eqref{eq-vp-psi-0}, we obtain an initial value problem for an integro-differential equation
\begin{equation}\label{eq-ivp-vp}
\left\{\!\begin{alignedat}{2}
& A_{11}(s,0)\vp'(s)+P(s,0)\vp(s)+\int_0^s\cL_{s,t}R(\si,0,s,0)\vp(\si)\,\rd\si=0, & \quad & -\ve\le s\le\ve,\\
& \vp(0)=0.
\end{alignedat}\right.
\end{equation}
By the unique existence of the solution (see, e.g., \cite{C98}), we have $\vp=0$ on $[-\ve,\ve]$. In a completely parallel manner, taking $s=0$ in \eqref{eq-ell-vp-psi} and using \eqref{eq-vp-psi-0} again yield
\begin{equation}\label{eq-ivp-psi}
\left\{\!\begin{alignedat}{2}
& A_{22}(0,t)\psi'(t)+Q(0,t)\psi(t)+\int_0^t\cL_{s,t}R(0,\tau,0,t)\psi(\tau)\,\rd\tau=0, & \quad & -\ve\le t\le\ve,\\
& \psi(0)=0,
\end{alignedat}\right.
\end{equation}
implying $\psi=0$ on $[-\ve,\ve]$. In view of formula \eqref{eq-sol-w}, we finally obtain $w=0$ on $\wt\om_1$ and thus
\[
u_2=0\quad\mbox{on }\om_1:=\{(x,y)\in\BR^2;\,(s(x,y),t(x,y))\in\wt\om_1\}.
\]
Since obviously $(x_0,y_0)\in\om_1\subset\om(x_0,y_0)$, we conclude $\bm u=(u_1,u_2)^\T=0$ in a smaller neighborhood $\om_1$ of $(x_0,y_0)$. Eventually, we finish the proof with the help of Lemma \ref{lem-NW}.

We close this section by illustrating the outline of the above proof in Figure 1. At the starting point, we have $u_1=0$ in $\om(x_0,y_0)$ and thus $u_2$ satisfies the overdetermined hyperbolic-elliptic system \eqref{eq-u2-hyper}--\eqref{eq-u2-ell} in $\om(x_0,y_0)$. Next, we perform the invertible changes of variables $s=s(x,y)$, $t=t(x,y)$ and introduce an auxiliary function $w(s(x,y),t(x,y)):=u_2(x,y)$, so that the hyperbolic equation that $w$ satisfies is reduced to \eqref{eq-w-hyper} on $\wt\om$. Employing the Riemann function of \eqref{eq-w-hyper} and the condition \eqref{eq-vanish-u2}, we can conclude the vanishing of $w$ on a subdomain $\wt\om_1\subset\wt\om$, which implies the vanishing of $u_2$ on a corresponding subdomain $\om_1\subset\om(x_0,y_0)$. Finally, the application of Lemma \ref{lem-NW} immediately completes the proof.\medskip

\centerline{%WinTpicVersion4.32a
{\unitlength 0.1in%
\begin{picture}(36.6900,16.5600)(0.0000,-16.5600)%
% STR 2 0 3 0 Black White  
% 4 610 334 610 354 1 0 0 0
% $(x_0,y_0)$
\put(6.1000,-3.5400){\makebox(0,0)[lt]{$(x_0,y_0)$}}%
% STR 2 0 3 0 Black White  
% 4 410 234 410 254 5 0 0 0
% $\omega$
\put(4.1000,-2.5400){\makebox(0,0){$\omega$}}%
% ELLIPSE 2 0 3 0 Black White  
% 4 610 354 110 54 110 54 110 54
% 
\special{pn 8}%
\special{ar 610 354 500 300 0.0000000 6.2831853}%
% DOT 0 0 3 0 Black White  
% 1 610 354
% 
\special{pn 4}%
\special{sh 1}%
\special{ar 610 354 16 16 0 6.2831853}%
% STR 2 0 3 0 Black White  
% 4 1810 234 1810 254 5 0 0 0
% $u_2(x,y)\to w(s,t)$
\put(18.1000,-2.5400){\makebox(0,0){$u_2(x,y)\to w(s,t)$}}%
% VECTOR 2 0 3 0 Black White  
% 2 1210 354 2410 354
% 
\special{pn 8}%
\special{pa 1210 354}%
\special{pa 2410 354}%
\special{fp}%
\special{sh 1}%
\special{pa 2410 354}%
\special{pa 2343 334}%
\special{pa 2357 354}%
\special{pa 2343 374}%
\special{pa 2410 354}%
\special{fp}%
% LINE 2 0 3 0 Black White  
% 4 2610 354 2860 104 3610 354 3360 102
% 
\special{pn 8}%
\special{pa 2610 354}%
\special{pa 2860 104}%
\special{fp}%
\special{pa 3610 354}%
\special{pa 3360 102}%
\special{fp}%
% CIRCLE 2 0 3 0 Black White  
% 4 3110 354 3360 104 3360 104 2860 104
% 
\special{pn 8}%
\special{ar 3110 354 354 354 3.9269908 5.4977871}%
% CIRCLE 2 0 3 0 Black White  
% 4 3110 854 3360 604 3360 604 2860 604
% 
\special{pn 8}%
\special{ar 3110 854 354 354 3.9269908 5.4977871}%
% CIRCLE 2 0 3 0 Black White  
% 4 2735 479 2610 354 2610 354 2860 604
% 
\special{pn 8}%
\special{ar 2735 479 177 177 0.7853982 3.9269908}%
% CIRCLE 2 0 3 0 Black White  
% 4 3485 479 3360 604 3360 604 3615 354
% 
\special{pn 8}%
\special{ar 3485 479 177 177 5.5173925 2.3561945}%
% DOT 0 0 3 0 Black White  
% 1 3110 254
% 
\special{pn 4}%
\special{sh 1}%
\special{ar 3110 254 16 16 0 6.2831853}%
% STR 2 0 3 0 Black White  
% 4 3110 234 3110 254 1 0 0 0
% $(0,0)$
\put(31.1000,-2.5400){\makebox(0,0)[lt]{$(0,0)$}}%
% STR 2 0 3 0 Black White  
% 4 2735 429 2735 479 5 0 0 0
% $\widetilde\omega$
\put(27.3500,-4.7900){\makebox(0,0){$\widetilde\omega$}}%
% STR 2 0 3 0 Black White  
% 4 267 1334 267 1354 5 0 0 0
% $\omega$
\put(2.6700,-13.5400){\makebox(0,0){$\omega$}}%
% ELLIPSE 2 0 3 0 Black White  
% 4 617 1354 117 1054 117 1054 117 1054
% 
\special{pn 8}%
\special{ar 617 1354 500 300 0.0000000 6.2831853}%
% DOT 0 0 3 0 Black White  
% 1 617 1354
% 
\special{pn 4}%
\special{sh 1}%
\special{ar 617 1354 16 16 0 6.2831853}%
% STR 2 0 3 0 Black White  
% 4 1817 1234 1817 1254 5 0 0 0
% $u_2(x,y)\leftarrow w(s,t)$
\put(18.1700,-12.5400){\makebox(0,0){$u_2(x,y)\leftarrow w(s,t)$}}%
% VECTOR 2 0 3 0 Black White  
% 2 2417 1354 1217 1354
% 
\special{pn 8}%
\special{pa 2417 1354}%
\special{pa 1217 1354}%
\special{fp}%
\special{sh 1}%
\special{pa 1217 1354}%
\special{pa 1284 1374}%
\special{pa 1270 1354}%
\special{pa 1284 1334}%
\special{pa 1217 1354}%
\special{fp}%
% LINE 2 0 3 0 Black White  
% 4 2617 1354 2867 1104 3617 1354 3367 1102
% 
\special{pn 8}%
\special{pa 2617 1354}%
\special{pa 2867 1104}%
\special{fp}%
\special{pa 3617 1354}%
\special{pa 3367 1102}%
\special{fp}%
% CIRCLE 2 0 3 0 Black White  
% 4 3117 1354 3367 1104 3367 1104 2867 1104
% 
\special{pn 8}%
\special{ar 3117 1354 354 354 3.9269908 5.4977871}%
% CIRCLE 2 0 3 0 Black White  
% 4 3117 1854 3367 1604 3367 1604 2867 1604
% 
\special{pn 8}%
\special{ar 3117 1854 354 354 3.9269908 5.4977871}%
% CIRCLE 2 0 3 0 Black White  
% 4 2742 1479 2617 1354 2617 1354 2867 1604
% 
\special{pn 8}%
\special{ar 2742 1479 177 177 0.7853982 3.9269908}%
% CIRCLE 2 0 3 0 Black White  
% 4 3492 1479 3367 1604 3367 1604 3622 1354
% 
\special{pn 8}%
\special{ar 3492 1479 177 177 5.5173925 2.3561945}%
% DOT 0 0 3 0 Black White  
% 1 3117 1254
% 
\special{pn 4}%
\special{sh 1}%
\special{ar 3117 1254 16 16 0 6.2831853}%
% STR 2 0 3 0 Black White  
% 4 2742 1429 2742 1479 5 0 0 0
% $\widetilde\omega$
\put(27.4200,-14.7900){\makebox(0,0){$\widetilde\omega$}}%
% BOX 2 0 3 0 Black White  
% 2 3292 1079 2942 1429
% 
\special{pn 8}%
\special{pa 3292 1079}%
\special{pa 2942 1079}%
\special{pa 2942 1429}%
\special{pa 3292 1429}%
\special{pa 3292 1079}%
\special{pa 2942 1079}%
\special{fp}%
% LINE 3 0 3 0 Black White  
% 38 3287 1204 3067 1424 3287 1168 3031 1424 3287 1132 2995 1424 3119 1264 2959 1424 3267 1080 2943 1404 3231 1080 2943 1368 3195 1080 2943 1332 3159 1080 2943 1296 3123 1080 2943 1260 3087 1080 2943 1224 3051 1080 2943 1188 3015 1080 2943 1152 2979 1080 2943 1116 3287 1096 3127 1256 3287 1240 3103 1424 3287 1276 3139 1424 3287 1312 3175 1424 3287 1348 3211 1424 3287 1384 3247 1424
% 
\special{pn 4}%
\special{pa 3287 1204}%
\special{pa 3067 1424}%
\special{fp}%
\special{pa 3287 1168}%
\special{pa 3031 1424}%
\special{fp}%
\special{pa 3287 1132}%
\special{pa 2995 1424}%
\special{fp}%
\special{pa 3119 1264}%
\special{pa 2959 1424}%
\special{fp}%
\special{pa 3267 1080}%
\special{pa 2943 1404}%
\special{fp}%
\special{pa 3231 1080}%
\special{pa 2943 1368}%
\special{fp}%
\special{pa 3195 1080}%
\special{pa 2943 1332}%
\special{fp}%
\special{pa 3159 1080}%
\special{pa 2943 1296}%
\special{fp}%
\special{pa 3123 1080}%
\special{pa 2943 1260}%
\special{fp}%
\special{pa 3087 1080}%
\special{pa 2943 1224}%
\special{fp}%
\special{pa 3051 1080}%
\special{pa 2943 1188}%
\special{fp}%
\special{pa 3015 1080}%
\special{pa 2943 1152}%
\special{fp}%
\special{pa 2979 1080}%
\special{pa 2943 1116}%
\special{fp}%
\special{pa 3287 1096}%
\special{pa 3127 1256}%
\special{fp}%
\special{pa 3287 1240}%
\special{pa 3103 1424}%
\special{fp}%
\special{pa 3287 1276}%
\special{pa 3139 1424}%
\special{fp}%
\special{pa 3287 1312}%
\special{pa 3175 1424}%
\special{fp}%
\special{pa 3287 1348}%
\special{pa 3211 1424}%
\special{fp}%
\special{pa 3287 1384}%
\special{pa 3247 1424}%
\special{fp}%
% STR 2 0 3 0 Black White  
% 4 3267 1084 3267 1104 4 0 0 0
% $\widetilde\omega_1$
\put(32.6700,-11.0400){\makebox(0,0)[rt]{$\widetilde\omega_1$}}%
% SPLINE 2 0 3 0 Black White  
% 3 339 1142 619 1212 893 1142
% 
\special{pn 8}%
\special{pa 339 1142}%
\special{pa 370 1154}%
\special{pa 432 1176}%
\special{pa 463 1186}%
\special{pa 494 1195}%
\special{pa 525 1202}%
\special{pa 556 1207}%
\special{pa 587 1211}%
\special{pa 618 1212}%
\special{pa 649 1211}%
\special{pa 680 1207}%
\special{pa 711 1201}%
\special{pa 742 1194}%
\special{pa 773 1185}%
\special{pa 804 1175}%
\special{pa 835 1164}%
\special{pa 866 1152}%
\special{pa 893 1142}%
\special{fp}%
% SPLINE 2 0 3 0 Black White  
% 3 893 1142 761 1362 739 1516
% 
\special{pn 8}%
\special{pa 893 1142}%
\special{pa 853 1194}%
\special{pa 834 1221}%
\special{pa 816 1248}%
\special{pa 799 1275}%
\special{pa 784 1303}%
\special{pa 772 1332}%
\special{pa 761 1361}%
\special{pa 753 1392}%
\special{pa 748 1423}%
\special{pa 744 1455}%
\special{pa 741 1487}%
\special{pa 739 1516}%
\special{fp}%
% SPLINE 2 0 3 0 Black White  
% 3 739 1516 617 1572 487 1516
% 
\special{pn 8}%
\special{pa 739 1516}%
\special{pa 710 1535}%
\special{pa 681 1552}%
\special{pa 652 1565}%
\special{pa 623 1572}%
\special{pa 593 1570}%
\special{pa 564 1561}%
\special{pa 535 1546}%
\special{pa 505 1528}%
\special{pa 487 1516}%
\special{fp}%
% SPLINE 2 0 3 0 Black White  
% 3 487 1516 469 1344 343 1144
% 
\special{pn 8}%
\special{pa 487 1516}%
\special{pa 486 1483}%
\special{pa 485 1451}%
\special{pa 483 1418}%
\special{pa 479 1387}%
\special{pa 472 1356}%
\special{pa 463 1326}%
\special{pa 451 1298}%
\special{pa 436 1270}%
\special{pa 419 1243}%
\special{pa 400 1216}%
\special{pa 360 1164}%
\special{pa 343 1144}%
\special{fp}%
% LINE 3 0 3 0 Black White  
% 42 783 1296 537 1542 761 1354 561 1554 749 1402 587 1564 741 1446 619 1568 737 1486 667 1556 845 1198 515 1528 843 1164 493 1514 609 1362 489 1482 735 1200 487 1448 691 1208 485 1414 649 1214 481 1382 613 1214 475 1352 577 1214 467 1324 545 1210 455 1300 515 1204 443 1276 485 1198 429 1254 459 1188 415 1232 431 1180 401 1210 405 1170 385 1190 379 1160 369 1170 787 1184 625 1346
% 
\special{pn 4}%
\special{pa 783 1296}%
\special{pa 537 1542}%
\special{fp}%
\special{pa 761 1354}%
\special{pa 561 1554}%
\special{fp}%
\special{pa 749 1402}%
\special{pa 587 1564}%
\special{fp}%
\special{pa 741 1446}%
\special{pa 619 1568}%
\special{fp}%
\special{pa 737 1486}%
\special{pa 667 1556}%
\special{fp}%
\special{pa 845 1198}%
\special{pa 515 1528}%
\special{fp}%
\special{pa 843 1164}%
\special{pa 493 1514}%
\special{fp}%
\special{pa 609 1362}%
\special{pa 489 1482}%
\special{fp}%
\special{pa 735 1200}%
\special{pa 487 1448}%
\special{fp}%
\special{pa 691 1208}%
\special{pa 485 1414}%
\special{fp}%
\special{pa 649 1214}%
\special{pa 481 1382}%
\special{fp}%
\special{pa 613 1214}%
\special{pa 475 1352}%
\special{fp}%
\special{pa 577 1214}%
\special{pa 467 1324}%
\special{fp}%
\special{pa 545 1210}%
\special{pa 455 1300}%
\special{fp}%
\special{pa 515 1204}%
\special{pa 443 1276}%
\special{fp}%
\special{pa 485 1198}%
\special{pa 429 1254}%
\special{fp}%
\special{pa 459 1188}%
\special{pa 415 1232}%
\special{fp}%
\special{pa 431 1180}%
\special{pa 401 1210}%
\special{fp}%
\special{pa 405 1170}%
\special{pa 385 1190}%
\special{fp}%
\special{pa 379 1160}%
\special{pa 369 1170}%
\special{fp}%
\special{pa 787 1184}%
\special{pa 625 1346}%
\special{fp}%
% STR 2 0 3 0 Black White  
% 4 617 1434 617 1454 5 0 0 0
% $\omega_1$
\put(6.1700,-14.5400){\makebox(0,0){$\omega_1$}}%
% VECTOR 2 0 3 0 Black White  
% 2 3100 700 3100 950
% 
\special{pn 8}%
\special{pa 3100 700}%
\special{pa 3100 950}%
\special{fp}%
\special{sh 1}%
\special{pa 3100 950}%
\special{pa 3120 883}%
\special{pa 3100 897}%
\special{pa 3080 883}%
\special{pa 3100 950}%
\special{fp}%
% STR 2 0 3 0 Black White  
% 4 3070 825 3070 875 3 0 0 0
% Riemann function
\put(30.7000,-8.7500){\makebox(0,0)[rb]{Riemann function}}%
\end{picture}}%
}\medskip
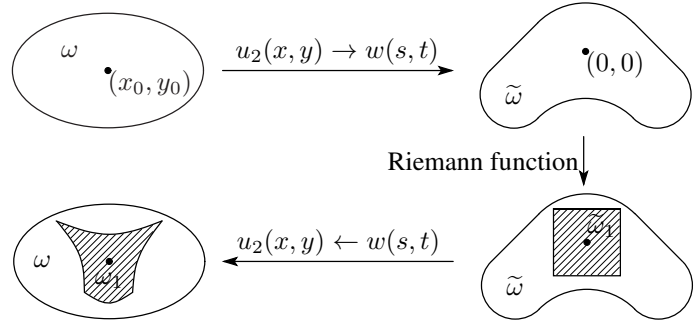
\centerline{\small {\bf Figure 1}.\quad Outline of the proof of Theorem \ref{thm-UCP}}

\section{Some Remarks and Examples}\label{sec-egrem}

In this section, we mainly construct various examples to illustrate the necessity of the additional information \eqref{eq-vanish-u2} as well as its possible reduction under special choices of coefficients.

Actually, by reviewing the argument of Step 2 in Section \ref{sec-proof}, it turns out that in principle we only need
\begin{equation}\label{eq-info-w-4}
w=\pa_s w=\pa_t w=0,\quad\pa_s^2w=0\mbox{ or }\pa_t^2w=0\quad\mbox{at }(0,0)
\end{equation}
instead of a stronger one \eqref{eq-vanish-w}. Indeed, without loss of generality, let us assume
\[
w=\pa_s w=\pa_t w=\pa_t^2w=0\quad\mbox{at }(0,0).
\]
Recalling the definition \eqref{eq-vp-psi} of $\vp$ and $\psi$, we find that the condition \eqref{eq-vp-psi-0} is now weakened as
\[
\vp(0)=\psi(0)=\psi'(0)=0.
\]
Even though, still we are able to derive the initial value problem \eqref{eq-ivp-vp} for $\vp$. This indicates $\vp=0$ on $[-\ve,\ve]$, and then simply taking $s=0$ in \eqref{eq-ell-vp-psi} again leads us to \eqref{eq-ivp-psi} and thus $\psi=0$ on $[-\ve,\ve]$.

On the other hand, it follows from \eqref{eq-w-hyper} that \eqref{eq-info-w-4} implies $\pa_s\pa_t w(0,0)=0$ immediately. Moreover, by the ellipticity of \eqref{eq-w-ell}, we see that \eqref{eq-info-w-4} automatically implies \eqref{eq-vanish-w}, that is, $w$ vanishes at $(0,0)$ up to all its second order derivatives. Therefore, it reveals that \eqref{eq-info-w-4} is equivalent to \eqref{eq-vanish-w} though the latter looks stronger. According to the relation \eqref{eq-u2-w-2}, this further results in
\begin{equation}\label{eq-info-u2-3}
\pa_x^2u_2=\pa_x\pa_y u_2=\pa_y^2u_2=0\quad\mbox{at }(x_0,y_0).
\end{equation}

Unfortunately, we will illustrate in the following example that the above facts no longer hold true on the opposite side. More precisely, it may happen that a similar assumption \eqref{eq-info-u2-4} of $u_2$ does not necessarily imply \eqref{eq-vanish-w}.

\begin{example}\label{ex-info-4}
We shall construct constant tensors $(a_{ijk\ell})$ satisfying conditions \eqref{eq-cond-ell} and \eqref{eq-cond-hyper} whereas \eqref{eq-info-u2-4} fails to imply \eqref{eq-vanish-w}.

(a)\ \ In the case of $u_2=\pa_x u_2=\pa_y u_2=\pa_x^2u_2=0$ at $(x_0,y_0)$, we choose
\[
a_{1111}\gg1,\quad a_{1112}=a_{1122}=0,\quad a_{1212}=2,\quad a_{1222}=a_{2222}=1.
\]
Then equations \eqref{eq-u2-hyper}--\eqref{eq-u2-ell} are linearly dependent at $(x_0,y_0)$ and we only have
\[
2\pa_x\pa_y u_2+\pa_y^2u_2=0\quad\mbox{at }(x_0,y_0),
\]
which allows nontrivial solutions.

(b)\ \ In the case of $u_2=\pa_x u_2=\pa_y u_2=\pa_y^2u_2=0$ at $(x_0,y_0)$, one possible choice can be
\[
a_{1111}\gg1,\quad a_{2222}\gg1,\quad a_{1212}=2,\quad a_{1122}=4,\quad a_{1112}=2,\quad a_{1222}=3.
\]
Then again \eqref{eq-u2-hyper}--\eqref{eq-u2-ell} give only one equation
\[
\pa_x^2u_2+3\pa_x\pa_y u_2=0\quad\mbox{at }(x_0,y_0).
\]

As long as two of the second derivatives of $u_2$ do not vanish at $(x_0,y_0)$, it is clear that the second derivatives of $w$ cannot vanish simultaneously at $(0,0)$ in view of \eqref{eq-u2-w-2}.
\end{example}

The above counterexamples result from the fact that such conditions on $(a_{ijk\ell})$ as \eqref{eq-cond-ell} and \eqref{eq-cond-hyper} are irrelevant of the linear dependency issue. Since it is unrealistic to make assumptions on $w(s,t)$ after changes of variables, we conclude from Example \ref{ex-info-4} that the original assumption \eqref{eq-vanish-u2} cannot be weakened e.g.\! to \eqref{eq-info-u2-4} in general.

Nevertheless, as was mentioned at the end of Section \ref{sec-main}, the vanishing assumption \eqref{eq-vanish-u2} can be reduced e.g.\! to \eqref{eq-info-u2-4} in the case of constant Lam\'e coefficients. Actually, the same reduction works as long as $a_{1222}(x_0,y_0)=0$ in the general case. To see this, we simply substitute \eqref{eq-info-u2-4} and $a_{1222}(x_0,y_0)=0$ into \eqref{eq-u2-hyper}--\eqref{eq-u2-ell} to deduce
\[
\begin{cases}
a_{1112}\pa_x^2u_2+(a_{1212}+a_{1122})\pa_x\pa_y u_2=0,\\
a_{1212}\pa_x^2+a_{2222}\pa_y^2u_2=0,
\end{cases}\mbox{at }(x_0,y_0).
\]
Since conditions \eqref{eq-cond-ell} and \eqref{eq-cond-hyper} guarantee
$
a_{1212}>0,  \ a_{2222}>0, \ a_{1212}+a_{1122}\ne0$ at $(x_0,y_0),$
it turns out that either $\pa_x^2u_2(x_0,y_0)=0$ or $\pa_y^2u_2(x_0,y_0)=0$ implies \eqref{eq-info-u2-3}, and thus all the subsequent arguments follow. As the condition $a_{1222}(x_0,y_0)=0$ definitely includes the orthotropic case, we see that the weakened assumption \eqref{eq-info-u2-4} works in a rather wide range of applications.

In this direction, we can even proceed further by specifying the coefficients in the governing equation \eqref{eq-gov-u}. In the following four examples, we will illustrate the possibility of determining $u_2$ by less than $4$ constants. We will restrict ourselves in the isotropic Lam\'e case \eqref{eq-Lame} possibly with some lower order terms, so that the assumptions in Lemma \ref{lem-NW} are fulfilled (see \cite[Example 2]{NW06}).

\begin{example}[$2$ constants]\label{ex-info-2}
In the isotropic Lam\'e system \eqref{eq-Lame}, we simply choose
$\mu(x,y)=\e^x, \ \la(x,y)=\e^y.$
Then $u_1=0$ in $\om(x_0,y_0)$ gives
\begin{equation}\label{eq-Lame-u2'}
\begin{cases}
(\e^x+\e^y)\pa_x\pa_y u_2=0,\\
\pa_x(\e^x\pa_x u_2)+\pa_y((2\,\e^x+\e^y)\pa_y u_2)=0
\end{cases}\mbox{in }\om(x_0,y_0).
\end{equation}
Obviously, the first equation of \eqref{eq-Lame-u2'} indicates that $u_2$ takes the form of separated variables
\begin{equation}\label{eq-u2-sep}
u_2(x,y)=p(x)+q(y).
\end{equation}
Plugging \eqref{eq-u2-sep} into the second equation of \eqref{eq-Lame-u2'} yields
\begin{equation}\label{eq-pq-2}
(\e^xp'(x))'+2\,\e^xq''(y)+(\e^yq'(y))'=0.
\end{equation}
Taking $y$-derivative in \eqref{eq-pq-2}, we obtain $2\,\e^xq'''(y)+(\e^yq'(y))''=0$. By the linear independence of $1$ and $\e^x$, we obtain $q'''(y)=(\e^yq'(y))''=0$. Then $q$ should be a quadratic function, i.e., $q(y)=C_0+C_1y+C_2y^2$. Further, we know
\[
0=(\e^yq'(y))''=\e^y((C_1+4C_2)+2C_2y)\quad\mbox{near }y_0
\]
and thus $C_1=C_2=0$. This means $q(y)$ is a constant, or equivalently $u_2(x,y)=p(x)$ only depends on $x$. Now \eqref{eq-pq-2} reduces to an ordinary differential equation with respect to $p$ near $x_0$, whose solution takes the form of
$
u_2(x,y)=p(x)=C_0+C_3\,\e^{-x}$ in $\om(x_0,y_0).$
In other words, $u_2$ only depends on $2$ constants $C_0,C_3$.
\end{example}

In the following three examples, we fix the two Lam\'e coefficients $\mu$ and $\la$ as constants satisfying $\mu+\la\ne0$, and set
\[
a_{1111}=a_{2222}=2\mu+\la,\quad a_{1212}=\mu,\quad a_{1122}=\la,\quad a_{1112}=a_{1222}=b_{121}=b_{122}=c_{12}=0.
\]
This means that we restrict $(a_{ijk\ell})$ in the constant Lam\'e case, but we add some nonzero lower order coefficients. Then $u_1=0$ in $\om(x_0,y_0)$ implies
\begin{align}
& (\mu+\la)\pa_x\pa_yu_2=0,\nonumber\\
& \mu\,\pa_x^2u_2+(2\mu+\la)\pa_y^2u_2+b_{221}\pa_xu_2+b_{222}\pa_yu_2+c_{22}u_2=0\label{eq-ell-u2}
\end{align}
in $\om(x_0,y_0)$. Again, it follows from the first equation above that $u_2$ takes the form \eqref{eq-u2-sep}. By adjusting $b_{221}$, $b_{222}$ and $c_{22}$, we explore other possibilities of determining $u_2$ by less constants than $4$.

\begin{example}[$3$ constants]
Choosing $b_{221}(x,y)=\e^y$ and $b_{222}=c_{22}=0$, we substitute \eqref{eq-u2-sep} into \eqref{eq-ell-u2} to deduce
\begin{equation}\label{eq-pq-3}
\mu\,p''(x)+(2\mu+\la)q''(y)+\e^yp'(x)=0\quad\mbox{in }\om(x_0,y_0).
\end{equation}
Differentiating with respect to $x$ gives $\mu\,p'''(x)+\e^yp''(x)=0$, and again the linear independence of $1$ and $\e^y$ implies $p''(x)=0$ or equivalently $p(x)=C_0+C_1x$. Similarly to Example \ref{ex-info-2}, we plug $p(x)$ back into \eqref{eq-pq-3} and arrive at an ordinary differential equation
\[
(2\mu+\la)q''(y)+C_1\e^y=0\quad\mbox{near} \ y_0.
\]
Therefore, it turns out that $u_2$ depends on $3$ constants $C_0,C_1,C_2$ in such a way that
\[
u_2(x,y)=C_0+C_1\left(x-\f{\e^y}{2\mu+\la}\right)+C_2y\quad\mbox{in }\om(x_0,y_0).
\]
\end{example}

In the next two examples, we set $x_0=y_0=0$ without loss of generality.

\begin{example}[$1$ constant]\label{ex-info-1}
Choosing $b_{221}(x,y)=xy$, $b_{222}(x,y)=xy^2$ and $c_{22}=0$, we have
\begin{equation}\label{eq-pq-1}
\mu\,p''(x)+(2\mu+\la)q''(y)+xy(p'(x)+y\,q'(y))=0\quad\mbox{in }\om(0,0).
\end{equation}
Taking $x=0$, \ $y=0$ and $x=y=0$ respectively in the above equation, we obtain
\begin{align*}
\mu\,p''(0)+(2\mu+\la)q''(y) & =0,\\
\mu\,p''(x)+(2\mu+\la)q''(0) & =0,\\
\mu\,p''(0)+(2\mu+\la)q''(0) & =0.
\end{align*}
We subtract the third identity from the summation of the first two to conclude
\begin{equation}\label{eq-pq-0}
\mu\,p''(x)+(2\mu+\la)q''(y)=0\quad\mbox{in }\om(0,0),
\end{equation}
which, in combination with \eqref{eq-pq-1}, yields $p'(x)+y\,q'(y)=0$ in $\om(0,0)$. Taking $y=0$ gives $p'(x)=0$ and further $q'(y)=0$. Consequently, both $p$ and $q$ are constants and so is $u_2$, i.e., $u_2(x,y)=C_0$.
\end{example}

Obviously, a constant is always a solution of \eqref{eq-ell-u2} if $c_{22}=0$. As a result, the last example achieves the unique determination of $u_2$ without any additional information by choosing a nonzero zeroth order term $c_{22}$.

\begin{example}[$0$ constant] \
Motivated by Example \ref{ex-info-1}, we simply choose $c_{22}(x,y)=xy$ and $b_{221}=b_{222}=0$, which gives
$
\mu\,p''(x)+(2\mu+\la)q''(y)+xy(p(x)+q(y))=0$ in $\om(0,0).$
Taking $x=0$, $y=0$ and $x=y=0$ respectively as before, again we can conclude \eqref{eq-pq-0}. However, this time we are directly led to
$
u_2(x,y)=p(x)+q(y)=0$ in $\om(0,0).$
In other words, the function $u_2$ does not depend on any constant and vanishes automatically.
\end{example}

\section{Conclusions}\label{sec-con}

In this paper, we established a new unique continuation property (UCP) for two-dimensional anisotropic elasticity systems in such a sense that the information of one component of the solution in a subdomain can uniquely determine the solution in the whole domain up to (at most) 5 constants. Our main Theorem \ref{thm-UCP} seems to be the first affirmative conclusion on this direction for anisotropic elasticity systems. Needless to say, such a result has remarkable practical significance because the required observation data can be almost halved in practice. Meanwhile, this also provides concrete evidence for the fact that small amounts of data sometimes result in satisfactory result in some practical problems.

As the first attempt to the UCPs with minimal information for elasticity systems, this paper is merely a starting point for a series of related problems awaiting investigation. Possible future topics include but by no means be restricted to the following aspects.
\begin{enumerate}
\item[(1)] For simplicity, throughout this paper we required $\bm u\in(C^2(\ov\Om))^2$ for the additional information \eqref{eq-vanish-u2} to make sense. However, it is achieved with rather smooth coefficients in view of the well-posedness of forward problems, which is usually invalid in such cases like Lam\'e coefficients with discontinuity. Hence, it is desirable to avoid the access to higher derivatives of $u_2$ at $(x_0,y_0)$ for practical applications. Therefore, an interesting problem is to find alternatives of \eqref{eq-vanish-u2}, e.g.\! $u_2=0$ at five distinct points.
\item[(2)] Compared with the two-dimensional case treated in this paper, the case of three-dimensional anisotropic elasticity systems for $\bm u=(u_1,u_2,u_3)^\T$ is definitely more challenging. We shall start from some special cases, e.g., the transversely isotropic case. Besides, the three-dimensional Riemann function also deserves closer studies.
\item[(3)] Since Theorem \ref{thm-UCP} is definitely a qualitative result, it is natural to further consider the quantitative UCP, that is, to estimate $\bm u$ in $\Om$ by some norm of $u_1$ in $\om(x_0,y_0)$ and $|\pa_\al u_2(x_0,y_0)|$ ($|\al|\le2$). If this is done, we can proceed to some related inverse and control problems as direct applications.
\item[(4)] In Section \ref{sec-egrem}, we constructed several examples to illustrate the possibility of further reducing the additional information \eqref{eq-vanish-u2}. In order to clarity the underlying mechanism, we shall look into details to find the conditions for the possible alleviation of additional information.
\item[(5)] In some applications, the surface Cauchy data $\{u,T u\}$ are taken as observation data, where $T$ is the stress operator on an accessible subboundary. It is important to establish similar unique continuation results by such boundary data, instead of the volume measurement of $u_1$ considered here, which should be a future topic.
\item[(6)] Finally, in the light of potential real-world applications, it is necessary to develop corresponding numerical methods and perform numerical verifications to support the theoretical achievement.
\end{enumerate}

\acknowledgements{The authors appreciate the valuable discussions with Gen Nakamura (Hokkaido University).
This work is supported by the A3 Foresight Program ``Modeling and Computation of Applied Inverse Problems'',
Japan Society for the Promotion of Science (JSPS) and National Natural Science Foundation of China (NSFC).
The first author is supported by NSFC (No.11971121).
The second and the fourth authors are partially supported by JSPS KAKENHI Grant Number JP15H05740.
The third and the fourth authors are supported by NSFC (No.11771270).
The fourth author is partly supported by NSFC (No.91730303) and RUDN University Program 5-100.
}

\bibliographystyle{unsrt}


\end{document}